\documentclass{svproc}
\usepackage{url}
\usepackage{amsmath,amssymb,mathrsfs}
\usepackage{graphicx}
\usepackage{booktabs}
\makeatletter
\renewcommand*\env@matrix[1][\arraystretch]{%
  \edef\arraystretch{#1}%
  \hskip -\arraycolsep
  \let\@ifnextchar\new@ifnextchar
  \array{*\c@MaxMatrixCols c}}
\makeatother

\begin{document}
\mainmatter              % start of a contribution
\title{Global Aerodynamic Design 
Optimization via Primal-Dual 
Aggregation Method}
\titlerunning{Primal-Dual 
Aggregation Method}  % abbreviated title (for running head)
\author{Emre \"{O}zkaya \and Nicolas R. Gauger}
\authorrunning{Emre \"{O}zkaya and Nicolas R. Gauger} % abbreviated author list (for running head)
\tocauthor{Emre \"{O}zkaya and Nicolas R. Gauger}
\institute{Chair for Scientific Computing, TU Kaiserslautern, \\Paul-Ehrlich-Strasse 34 67663 Kaiserslautern, Germany\\
\email{\{emre.oezkaya,nicolas.gauger\}@scicomp.uni-kl.de},\\ WWW home page:
\texttt{https://www.scicomp.uni-kl.de}}

\maketitle              % typeset the title of the contribution

\begin{abstract}
Global aerodynamic design optimization using Euler or Navier-Stokes equations requires very reliable surrogate modeling techniques since the computational effort for the underlying flow simulations is usually high. In general, for such problems, the number of samples that can be generated to train the surrogate models is very limited due to restricted computational resources. On the other hand, recent developments in adjoint methods enable nowadays evaluation of gradient information at a reasonable computational cost for a wide variety of engineering problems. Therefore, a much richer data set can be collected using an adjoint solver in a Design of Experiment framework. In the present work, we present a novel aggregation method, which enables the state of the art surrogate models to incorporate extra gradient information without causing overfitting problems. Therefore, accurate surrogate models with relatively large number of design parameters can be trained from a small set of samples. We also present results two well known benchmark design optimization problems showing efficiency and robustness of the new method.
\end{abstract}
\section{Introduction}

For many engineering design problems, there is a growing tendency in the industry to design products using simulation based optimization.  As good examples of application fields, aerodynamic shape optimization and optimal flow control can be given. Typically, Euler, Navier-Stokes or RANS (Reynolds-Averaged Navier-Stokes) type of flow solvers are utilized to simulate the flow field and predict the aerodynamic performance. Since these kind of simulations may be computationally very expensive (e.g., several hundreds of CPU hours for a 3-$D$ simulation of a car body), surrogate modeling has become an important topic since it facilitates a very fast evaluation of a design configuration. Therefore, employing global optimization methods that may require many design evaluations becomes feasible.   

Especially in the last two decades, developments in adjoint methods has
made the efficient evaluation of gradient vectors feasible. For example, by using continuous or discrete adjoint flow solvers, it is possible to compute the full gradient vector of a given objective function within a computational effort that is only a small multiple of the computational effort required for the flow simulation itself. Therefore, a rich data set can be generated in the sampling stage since the adjoint solver can evaluate many design sensitivities as well as the objective itself. Therefore, it is desired to extend the existing surrogate modeling techniques so that they can benefit from the gradient information to improve the model accuracy. In this way, much better performance in the optimization stage can be achieved for the same computational budget.    

In this article, our main focus is a new surrogate modeling technique referred to as the primal-dual aggregation model. First, we briefly introduce the state of the art surrogate modeling techniques that use primal and dual information. Then, we present the new model and discuss its advantages over the present models. Finally, we present results of aerodynamic shape optimization studies for two benchmark test cases to show feasibility and efficiency of the suggested method.  

\section{Surrogate Modeling using Primal and Dual Information}
\label{sec:1}
In the most general setting, we have some data obtained from expensive flow simulations performed at different design configurations. In general, this initial data acquisition phase can be realized by using any proper Design of Experiment (DoE) method such as the Latin Hypercube Sampling (LHS) method. Note that, with the increasing number of design parameters, it becomes more and more difficult to realize a "space filling" DoE plan. Nevertheless, to begin with, we assume that the  functional values (primal information) at $N$ different design points are available: 
\begin{equation}
y^{(i)} = f(\bf x^{(i)}), \; i=1,\ldots, N.
\label{eq:yi}
\end{equation} 
Furthermore, we also assume that at the same design points we have also the gradient information (dual information) obtained by an adjoint solver:
\begin{equation}
\nabla y^{(i)} = \bar f(\bf x^{(i)}), \; i=1,\ldots, \overline N,
\label{eq:del_yi}
\end{equation}
where $f(\bf x)$ and $\bar f(\bf x)$ denote primal and adjoint simulations performed at a certain design configuration $\bf x$. Ideally, adjoint simulations are performed exactly at all design points specified by the DoE method. In this case, we have $\overline N = N$ and the training data can be obtained using only the adjoint solver. In a practical situation, however, the number of samples with gradient information may be less since the adjoint solution may not be available at every design point ($\overline N < N$). This is primarily due to the fact that, adjoint solvers are generally much more susceptible to convergence problems than the primal non-linear solvers. Therefore, the adjoint solution may not be available at each design point in a realistic optimization study. 

In the following, we assume that the target function $f$ is modeled by the well known Kriging model (also known the Gaussian process regression model). Note that, the aggregation model does not rely on a particular model and it can be constructed using any surrogate model trained by the primal functional values $y^{(i)}$. In the Kriging model, the output $f$ assumed to be a sum of a global trend function $G(\mathbf{x})$ and an error term $\epsilon(\mathbf{x})$:
\begin{equation}
\tilde f(\mathbf{x}) = G(\mathbf{x}) + \epsilon(\mathbf{x}).
\end{equation}
The key assumption in the Kriging model is that the error term $\epsilon(\mathbf{x})$ is a Gaussian process and the error at a point is correlated with the error at other points. For regression tasks, for which the training data comes from deterministic simulations, this assumption makes perfectly sense. As far as the global trend function is concerned, ideally a constant or a first order linear term can be used. If the global trend function is a constant term ($G(\mathbf{x}) = \beta_0$), we have then the so-called ``ordinary Kriging" model as the primal surrogate model.

The error term $\epsilon(\mathbf{x})$ is assumed to be a stationary random process with a zero mean, positive variance and covariance. As we mentioned previously, the error terms are spatially correlated with each other in the design space. The exact correlation, however, is unknown. To bridge the gap, the spatial correlation is assumed to be represented by certain analytic correlation functions. The most frequently used correlation functions are the radial basis functions:
\begin{equation}
R(x^{(i)},x^{(j)}) = e^{-\sum_{k=1}^d \theta_k | x^{(i)}_k-x^{(j)}_k |^{\gamma_k}},
\end{equation}
where $x^{(i)}_k$ and $x^{(j)}_k$ denote the $k$th elements of the $d$ dimensional design vectors $x^{(i)}$ and $x^{(j)}$ respectively. The parameters of the correlation function $\theta_k > 0$ and $0 < \gamma_k <2$ for $k=1,\ldots, d$ must be tuned to find the setting, which provides the maximum likelihood of the Kriging estimator. Note that as the distance between two designs $x^{(i)}$ and $x^{(j)}$ grows, the correlation function goes to zero so that the correlation between two points is very weak. In the other extreme case, the correlation function approaches to the value of one as the distance between $x^{(i)}$ and $x^{(j)}$ gets smaller.

The Kriging model assumes that the output at the sampled points $y(\bf x^{(i)}), \; i=1,\ldots N$ are drawn from a Gaussian process with a joint probability density function:
\begin{equation}
pdf(\mu,\Sigma)=\frac{1}{\sqrt{(2\pi)^N det(\Sigma)}} e^{-\frac{1}{2} (\mathbf{y} - \bf \mu)^\top \Sigma^{-1} (\mathbf{y}- \bf \mu) }.
\end{equation}
Here $\mathbf{\mu}$ is the $N$ dimensional mean vector $\mu = \mathbb{E}[\mathbf{y}]= [\mathbb{E}[y^{(1)}],\mathbb{E}[y^{(2)}], \ldots, \mathbb{E}[y^{(N)}] ] = [\beta_0, \beta_0, \ldots, \beta_0]$. In the following, this vector is represented as $\mathbf{1} \beta_0$, where $\mathbf{1}$ is a vector full of ones. The $\Sigma$ is on the other hand is the $N \times N$ covariance matrix.

If we add a ``new point" $(\mathbf{x^{(\ast)}}, y^{(\ast)})$ into the data set and we assume that the extended set
\begin{equation}
\mathbf{y_{new}} = [y^{(1)},y^{(2)}, \ldots ,y^{(N)}, y^{(\ast)} ]
\end{equation}
is also a $N+1$ variate Gaussian distribution, the mean value of the new distribution is simply the Kriging predictor we look for. It is given as
\begin{equation}
\beta_0^{\ast} = \beta_0 + \mathbf{\Sigma^{\ast}}^\top \Sigma^{-1} (\mathbf{y} - \mathbf{1} \beta_0), 
\label{eq:Kriging_cov}
\end{equation}
where $\mathbf{\Sigma^{\ast}}$ is the vector of covariances between the new point $y^{(\ast)}$ and the vector $\mathbf{y}$. Since we assume each element of $\mathbf{y}$ has the same variance (i.e., $\mathbf{y}$ is identically distributed), we can use correlations in Eq. \eqref{eq:Kriging_cov} instead of covariances:
\begin{equation}
\beta_0^{\ast} = \beta_0 + \mathbf{r}^\top \mathcal{R}^{-1} (\mathbf{y} - \mathbf{1} \beta_0),
\label{eq:Kriging_cor}
\end{equation}
where the correlation matrix $\mathcal{R}$ is given as 
\begin{equation}
\mathcal{R}= 
\begin{bmatrix}
R(\mathbf{x^{(1)}},\mathbf{x^{(1)}}) & R(\mathbf{x^{(1)}},\mathbf{x^{(2)}}) & \ldots & R(\mathbf{x^{(1)}},\mathbf{x^{(N)}}) \\
R(\mathbf{x^{(2)}},\mathbf{x^{(1)}}) & R(\mathbf{x^{(2)}},\mathbf{x^{(2)}}) & \ldots & R(\mathbf{x^{(2)}},\mathbf{x^{(N)}}) \\
\vdots & \vdots & \vdots & \vdots \\
R(\mathbf{x^{(N)}},\mathbf{x^{(1)}}) & R(\mathbf{x^{(N)}},\mathbf{x^{(2)}}) & \ldots & R(\mathbf{x^{(N)}},\mathbf{x^{(N)}}) 
\end{bmatrix}.
\end{equation}
Furthermore, $\mathbf{r}$ is the vector of correlations between the new point $y^{(\ast)}$ and the vector $\mathbf{y}$:
\begin{equation}
\mathbf{r} = [R(\mathbf{x^{(1)}},\mathbf{x^{(\ast)}}), R(\mathbf{x^{(2)}},\mathbf{x^{(\ast)}}), \ldots, R(\mathbf{x^{(N)}},\mathbf{x^{(\ast)}})].
\end{equation}
Since the $R(\mathbf{x^{(i)}},\mathbf{x^{(j)}})= R(\mathbf{x^{(j)}} ,\mathbf{x^{(i)}}) > 0, \; 0 \leq i,j \leq N$, the correlation matrix 
$\mathcal{R}$ is a symmetric positive definite matrix. All its diagonal entries are one and all its off diagonal entries are between zero and one. This property is very useful since efficient linear algebra algorithms such as the Cholesky decomposition can be used to solve the linear systems $\mathcal{R} w = (\mathbf{y} - \mathbf{1} \beta_0),\; w \in \mathbb{R}^N$ to evaluate the Kriging estimator given in Eq. \eqref{eq:Kriging_cor}. The bias term $\beta_0$, on the other hand is found by taking the least squares estimate, for which the details are skipped here for brevity. 

The Kriging model can be further enhanced to incorporate the gradient information. In the so-called indirect Gradient Enhanced Kriging (GEK) approach, new samples are generated using the first order Taylor approximation around each sample, where the functional and gradient information are already known \cite{Indirect_GEK}. Using this augmented data set, the Kriging model is trained as usual. Although the its implementation is quite straightforward, the indirect GEK method has some serious drawbacks. Since the Taylor approximation produces reliable results only in the close neighborhood of a sample, only a small perturbation in each direction can be used and the new samples tend to cluster around the initial samples. This worsens the condition number of the correlation matrix heavily and may lead to serious robustness problems. Another disadvantage of the indirect GEK model is that, the size of the correlation matrix gets rapidly larger with the increasing number of design parameters and number of samples, rendering the model training a very challenging problem. Probably due these shortcomings, the indirect GEK approach did not gain much popularity in the past. 

The most popular approach to incorporate gradient information to the Kriging model is the so-called direct GEK approach \cite{Han_VFM,Laurenceau}, in which the correlation matrix is augmented using the correlations of functional and sensitivity values with each other. Thereby, the main assumption is that the covariance between the functional values and the its derivatives are the derivative of the covariance with respect to the same variable:
\begin{equation}
cov \left( y(\mathbf{x^{(i)}}), \frac{\partial y(\mathbf{x^{(j)}} )}{\partial x_k^{(j)}} \right) = \frac{\partial}{\partial x_k^{(j)}} \left( cov(y(\mathbf{x^{(i)}}),y(\mathbf{x^{(j)}} ) \right).
\end{equation}

%The extended correlation matrix in the direct GEK method is given as
%\begin{equation}
%\dot{\mathcal{R}} = 
%\begin{bmatrix}[1.5]
%\mathcal{R} & \frac{\partial \mathcal{R}}{\partial x_1^{(i)}} & \frac{\partial \mathcal{R}}{\partial x_2^{(i)}} & \ldots &  \frac{\partial \mathcal{R}}{\partial x_d^{(i)}} \\
%\frac{\partial \mathcal{R}}{\partial x_1^{(j)}} & \frac{\partial^2 \mathcal{R}}{\partial x_1^{(i)}\partial x_1^{(j)}} & \frac{\partial^2 \mathcal{R}}{\partial x_1^{(i)}\partial x_2^{(j)}} & \ldots & \frac{\partial^2 \mathcal{R}}{\partial x_1^{(i)}\partial x_d^{(j)}} \\
%\frac{\partial \mathcal{R}}{\partial x_2^{(j)}} & \frac{\partial^2 \mathcal{R}}{\partial x_1^{(i)}\partial x_2^{(j)}} & \frac{\partial^2 \mathcal{R}}{\partial x_2^{(i)}\partial x_2^{(j)}} & \ldots & \frac{\partial^2 \mathcal{R}}{\partial x_2^{(i)}\partial x_d^{(j)}} \\
%\vdots \\  
%\frac{\partial \mathcal{R}}{\partial x_d^{(j)}} & \frac{\partial^2 \mathcal{R}}{\partial x_1^{(i)}\partial x_d^{(j)}} & \frac{\partial^2 \mathcal{R}}{\partial x_2^{(i)}\partial x_d^{(j)}} & \ldots & \frac{\partial^2 \mathcal{R}}{\partial x_d^{(i)}\partial x_d^{(j)}}
%\end{bmatrix}
%\label{eq:R_GEK}
%\end{equation}

The direct GEK approach requires correlation functions, which are at least two times differentiable. Therefore, the set of correlation function that can be used in this model is rather restricted. Commonly, the Gaussian correlation function is used ($\gamma = 2$). One can obtain with the direct GEK very good results and improve significantly the accuracy of the surrogate model with the same number of samples. It has, however, two major drawbacks:
\begin{itemize}
\item Similar to the indirect GEK method, the correlation matrix $\mathcal{R}$ grows rapidly as the number of samples ($N$) and design parameters ($d$) increase. Thus, direct GEK is not a scalable method. For example, with a moderate problem size of $N=500$ and $d=50$, the correlation matrix is a $25500 \times 25500$ dense matrix that requires approximately $5$Gb space in memory. Consequently, training the model hyper-parameters becomes computationally unfeasible. This issue has been addressed by Han et al. in \cite{WGEK}. The authors suggested training the direct GEK model only with one gradient vector at a time and taking a weighted average thorough all samples.

\item The accuracy of the GEK model is very vulnerable to the noise in the data. For smooth data, the direct GEK is capable of producing estimations at a very high accuracy. With the increasing noise in data, however, the generalization error grows rapidly leading to situations that direct GEK may generalize even worse than the Kriging model. This is due to the fact that high frequency components of the target function generates very high gradients that lead to very ill-conditioned correlation matrices. With regularization methods, such as adding a small multiple of the identity matrix to the correlation matrix, this problem can be partially circumvented. Finding the correct regularization, however, can be very challenging as there is usually no prior knowledge about the amount of noise in the data.     
\end{itemize}

\section{The Primal-Dual Aggregation Model}
The aggregation model, which is in essence a type of boosting method, is a convex combination of the primal and dual models with a scalar weight parameter $\alpha$:
\begin{equation}
\tilde f_{agg} (\mathbf{x}) = \alpha \tilde f_{dual} (\mathbf{x})+ (1-\alpha)\tilde f_{primal} (\mathbf{x}), \; \alpha \in [0,1]. 
\label{eq:primal_dual_agg}
\end{equation} 

As the primal model, we simply take the Kriging approach since it is proven to have a low generalization (out of sample) error for small data sets. For the dual model, we use the Taylor approximation around the nearest sample point, which is found by the nearest neighbor search algorithm:
\begin{equation}
\tilde f_{dual} (\mathbf{x}) = f(\mathbf{x_g}) + \nabla f(\mathbf{x_g})^\top (\mathbf{x}-\mathbf{x_g})
\end{equation}

The weight parameter $\alpha$ in the Eq. \eqref{eq:primal_dual_agg} is a spatially varying parameter (in design space) and it is a function of $L^1$ distance to the nearest point, where gradient information is available ($\mathbf{x_g}$), and the $L^1$-norm of the gradient at that point:
\begin{equation}
\alpha(\mathbf{x}) = e^{-\rho \|\mathbf{x} - \mathbf{x_g} \|_1 \|\nabla f(\mathbf{x_g}) \|_1 },
\end{equation}
where $\rho \in \mathbb{R}^+$ is a global hyper-parameter of the aggregation model. The most suitable value for this parameter is found by training the aggregation model using the K-fold cross validation technique. As the distance between $\mathbf{x}$ and $\mathbf{x_g}$ increases, $\alpha$ decreases. This means, the aggregation model produces essentially same results as the primal model if the gradient information is available only at a very distant sample point. On the other extreme, as $\mathbf{x}$ approaches $\mathbf{x_g}$, $\alpha$ tends to go to one. In this case, the aggregation model becomes exactly the first order Taylor approximation. In regions in the design space, where there is a sample point nearby, the dual model tends to have a larger weight. Similar to the distance $\|\mathbf{x} - \mathbf{x_g} \|_1$, the $L1$ norm of the gradient vector plays a similar role. The $L^1$ norm is chosen since it has a higher variation that leads to a better contrast with high dimensional data. The aggregation model has two distinct advantages compared to the GEK model:
\begin{itemize}
\item The aggregation method is much less susceptible to noise or inaccuracies in the gradient computation compared to the GEK method. This is particularly important for global optimization in aerodynamics since it is almost unavoidable that adjoint solvers generate noisy or inaccurate solutions at some design configurations. Therefore, the generalization error is lower, which leads to an optimization process that is more robust and reliable.
\item Higher number of design parameters and samples can be used since the hyper-parameter training of the aggregation model requires much less computational effort than the GEK. This is due to the fact that size of the correlation matrix increases very moderately with the increasing $N$ and $d$.
\end{itemize}

\section{Aerodynamic Shape Optimization Results}

In order to demonstrate the feasibility and efficiency of the aggregation model, shape optimization studies on generic configurations are performed. The primal and the adjoint solver used in the computations is the open source multi-physics suite SU$^2$ \cite{SU2}, which was originally developed at the Stanford University. The SU$^2$ code is based on the Finite-Volume method and can solve compressible Euler/RANS equations. For the evaluation of gradients, both continuous and discrete adjoint versions are available. In the present work, only the discrete adjoint version \cite{SU2_disc_adj} has been employed. In the optimization process, the initial DoE stage is realized by the LHS method. The actual optimization is performed using the Efficient Global Optimization (EGO) approach \cite{DACE} using the Expected Improvement (EI) function \cite{EI_method} to select the most promising designs. The optimization constraints are treated by training extra surrogate models for each of them and penalizing unfeasible designs with large penalty coefficients. The stopping criteria for the optimization is a predefined computational budget, which is specified as the maximum number of simulations allowed in total.     

\subsection{Transonic Drag Minimization}
\begin{figure}[h!]
  \centering
 \includegraphics[scale=0.65]{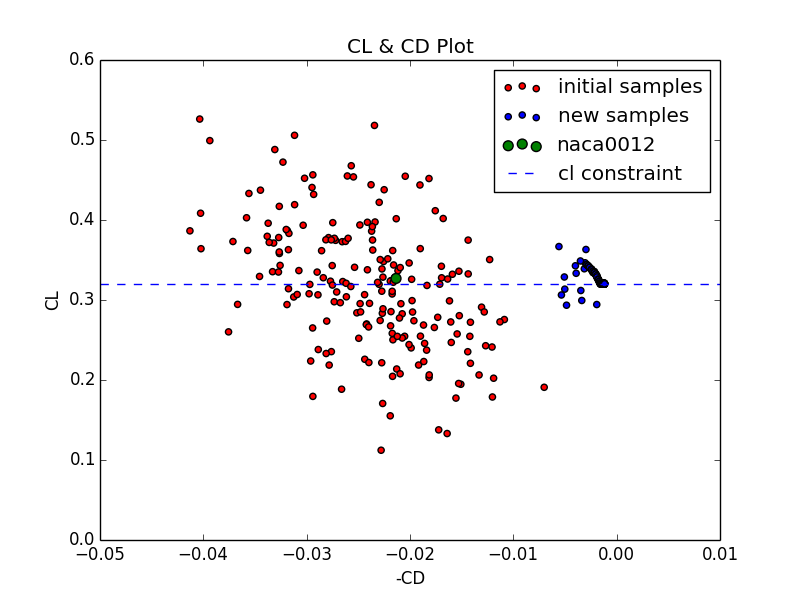} 
 \caption{NACA0012 drag minimization with lift constraint} 
 \label{fig:naca0012}
\end{figure}

As the first test case, we consider the drag minimization of NACA0012 airfoil at transonic flow conditions ($Ma=0.85$) with an Angle of Attack (AoA) of $2^\circ$. The computations are carried out using an unstructured grid with $10216$ control volumes. The shape parameters are defined to be the coefficients of $38$ Hicks-Henne basis functions, which are applied to deform upper and lower surfaces of the airfoil. For the simulations, $2$D Euler solver of the SU$^2$ with the Jameson-Schmidt-Turkel (JST) scheme is used. The lift coefficient and the airfoil area of the NACA0012 airfoil subtracted by a small tolerance are given as the inequality constraint to the problem such that the designs having smaller lift coefficient or/and area are marked as unfeasible designs. In Fig. \ref{fig:naca0012}, the $C_L$ and $C_D$ values of the samples that are generated in the optimization process are plotted. It can be observed that, the drag can be significantly reduced already after the initial data acquisition stage (red samples). The global optimization process using the primal-dual aggregation method improves the shape further (blue samples). We have also made a comparison between the optimization result obtained by the Kriging method and the suggested new method. The $C_L$ and $C_D$ values of the best design as well as the result obtained from the standard Kriging method (without gradients) are tabulated in Table \ref{tbl:naca0012}. It can be seen that, significantly better results are obtained by using objective and constraint gradients in the optimization process.

\begin{table}
 \begin{center}
  \begin{tabular}{lcccccc}
\hline 
    \textbf{configuration} & \textbf{Lift}    & \textbf{Drag}  & \textbf{Area}  \\ 
 \hline  \\
    NACA0012       & $0.3269\;$                 & $0.02133\;$    &    $0.08169$                       \\
    optimal (aggregation)       & $0.3201\;$                 & $0.00117\;$    &    $0.08104$                                         \\
    optimal (Kriging)    & $0.3301$               & $0.00283$    &    $0.08194$                              \\\\
  \end{tabular}
  \caption{\label{tbl:naca0012} Minimization results obtained for the NACA0012 profile.}
 \end{center}
\end{table}

\subsection{Viscous Drag Minimization}
\begin{figure}[h!]
 \centering
 \includegraphics[scale=0.65]{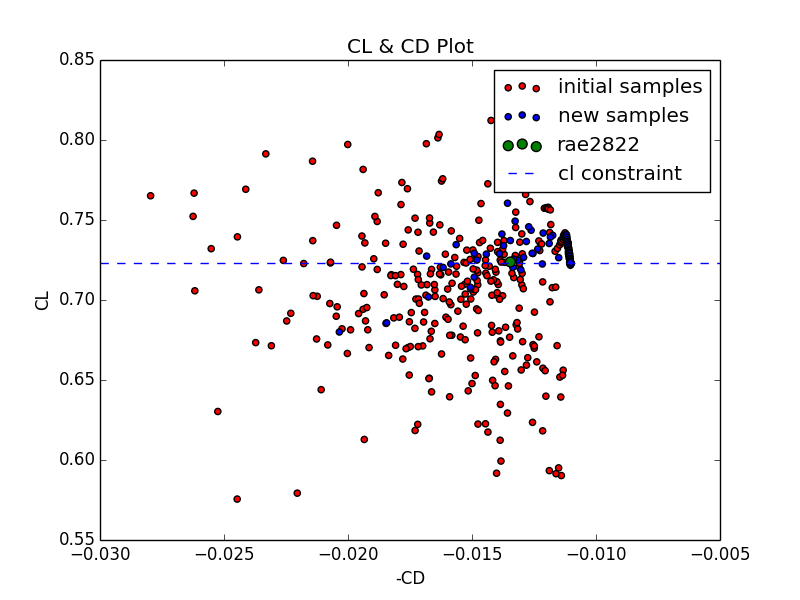} 
 \caption{RAE2822 drag minimization with lift constraint} 
 \label{fig:rae2822}
\end{figure}

The second test case is the drag minimization problem of the RAE2822 airfoil at $Ma=0.729$, $Re=6.5E06$ and AoA=$2.31^\circ$. The simulations are carried out solving $2$D RANS equations with the Spalart-Allmaras turbulence model. For the spatial discretization JST scheme is used. The computational domain is a hybrid mesh consisting of a structured part for the boundary layer and an unstructured part for the outer layer with total $22842$ control volumes. The shape parameters are coefficients of the $38$ Hicks-Henne basis functions. The lift and area constraints are set in a similar way to the previous test case. In Fig. \ref{fig:rae2822}, the $C_L$ and $C_D$ values of the samples that are generated in the optimization process are plotted. In contrast to the NACA0012 test case, the initial data acquisition stage did not bring a significant improvement. The $C_L$ and $C_D$ values of the best design as well as the result obtained from the standard Kriging method (without gradients) are tabulated in Table \ref{tbl:rae2822}. Similar to the previous test case, better results compared to Kriging are obtained by using the aggregation method using gradient information.

\begin{table}
 \begin{center}
  \begin{tabular}{lcccccc}
\hline 
    \textbf{configuration} & \textbf{Lift}    & \textbf{Drag}  & \textbf{Area}  \\ 
 \hline  \\
    RAE2822       & $0.7237\;$                 & $0.01345\;$    &    $0.07781$                       \\
    optimal (aggregation)    & $0.7222\;$               & $0.01101\;$    &    $0.07780$                              \\
    optimal (Kriging)    & $0.7275\;$               & $0.01182\;$    &    $0.07827$                              \\\\
  \end{tabular}
  \caption{\label{tbl:rae2822} Minimization results obtained for the RAE2822 profile.}
 \end{center}
\end{table}

\section{Conclusions}
We presented a novel primal-dual aggregation method for the global aerodynamic design optimization using expensive CFD simulations. The suggested method is a convex combination of two surrogate models with a spatially varying weight parameter that is tuned by the cross-validation technique. The main advantages of the suggested method are its robustness and low training effort, which makes it suitable for design optimization studies in industry scale. We have tested the aggregation method using two benchmark shape optimization problems of the SU$^2$ test case repository. It has been observed that, in both cases, the optimization process ran smoothly leading to significantly improved designs compared to the classical Kriging method.


\end{document}